\documentclass{gtart}

\def\ifplaintex{\expandafter\ifx\csname documentclass\endcsname\relax}

\def\gtp{{\mathsurround=0pt\it $\cal G\mskip-2mu$eometry \&\ 
$\cal T\!\!$opology $\cal P\!$ublications}}  

\def\Addressesr{\bigskip
{\small \parskip 0pt \leftskip 0pt \rightskip 0pt plus 1fil \def\\{\par}
\sl\theaddress\par
\medskip
\rm Email:\stdspace\tt\theemail\hfill\rm Received:\qua\receiveddate \par}}

\def\recd{{\small Received:\qua\receiveddate\ifx\reviseddate\relax
\else\qquad Revised:\qua\reviseddate\fi\par}} 

\def\lognumber#1{\def\thelognumber{#1}}
\def\volumenumber#1{\def\thevolumenumber{#1}}
\def\volumeyear#1{\def\thevolumeyear{#1}}
\def\papernumber#1{\def\thepapernumber{#1}}
\def\pagenumbers#1#2{\def\startpage{#1}\def\finishpage{#2}}
\def\published#1{\def\publishdate{#1}}

\def\received#1{\def\receiveddate{#1}}

\def\accepted#1{\def\accepteddate{#1}}

\let\\\par\let\thelognumber\relax\let\thevolumenumber\relax
\let\thepapernumber\relax\let\thevolumeyear\relax\let\startpage\relax
\let\finishpage\relax\let\publishdate\relax\let\receiveddate\relax
\let\reviseddate\relax\let\accepteddate\relax\let\theasciititle\relax
\let\theasciiauthors\relax
\let\theasciiabstract\relax

\let\theasciiemail\relax

\ifplaintex
\font\logobig=cmssbx10 scaled 3836
\font\logomed=cmssbx10 scaled 2557
\else
\font\logobig=cmssbx10 scaled 4200
\font\logomed=cmssbx10 scaled 2800
\fi

\long\def\makeagttitle{   
\count0=\startpage
\agt\hfill      
\hbox to 45truept{\vbox to 0pt{\vglue -13truept{\logomed A\kern -.37em{\logobig 
T}\kern -.38em G}\vss}\hss}
\break
{\small Volume \thevolumenumber\ (\thevolumeyear)
\startpage--\finishpage\nl
Published: \publishdate}

\vglue .25truein

{\parskip=0pt\leftskip 0pt plus
1fil\def\\{\par\smallskip}{\Large\bf\thetitle}\par\medskip} \vglue
0.05truein

{\parskip=0pt\leftskip 0pt plus 1fil\def\\{\par}{\sc\theauthors}
\par\medskip}%
 
\vglue 0.03truein 

{\small\leftskip 25truept\rightskip 25truept{\bf Abstract}\stdspace\theabstract

{\bf AMS Classification}\stdspace\theprimaryclass
\ifx\thesecondaryclass\relax\else; \thesecondaryclass\fi\par
{\bf Keywords}\stdspace \thekeywords\par}\vglue 7truept

}   

\ifplaintex
\font\phead=cmsl9 scaled 950
\font\pnum=cmbx10 scaled 913
\font\pfoot=cmsl9 scaled 950
\headline{\vbox to 0pt{\vskip -4.5mm\line{\small\phead\ifnum
\count0=\startpage ISSN 1472-2739 (on-line) 1472-2747 (printed)
\hfill {\pnum\folio}\else\ifodd\count0\def\\{ }%
\ifx\theshorttitle\relax\thetitle\else\theshorttitle\fi\hfill{\pnum\folio}
\else\def\\{ and }{\pnum\folio}\hfill\ifx\theshortauthors\relax\theauthors
\else\theshortauthors\fi\fi\fi}\vss}}
\footline{\vbox to 0pt{\vglue 0mm\line{\small\pfoot\ifnum\count0=\startpage
\copyright\ \gtp\hfill\else
\agt, Volume \thevolumenumber\ (\thevolumeyear)\hfill\fi}\vss}}
\else
\font=cmsl9 scaled 1050
\font\lnum=cmbx10 
\font=cmsl9 scaled 1050
\makeatletter
\def\@oddhead{{\small\lhead\ifnum\count0=\startpage ISSN 1472-2739 
(on-line) 1472-2747 (printed)\hfill {\lnum\number\count0}\else\ifodd\count0
\def\\{ }\ifx\theshorttitle\relax \thetitle \else\theshorttitle\fi\hfill
{\lnum\number\count0}\else\def\\{ and }{\lnum\number\count0}
\hfill\ifx\theshortauthors\relax 
\theauthors\else\theshortauthors\fi\fi\fi}}\def\@evenhead{\@oddhead}
\def\@oddfoot{\small\lfoot\ifnum\count0=\startpage\copyright\ \gtp\hfill\else
\agt, Volume \thevolumenumber\ (\thevolumeyear)\hfill\fi}
\def\@evenfoot{\@oddfoot}
\makeatother
\fi
\let\maketitlepage\makeagttitle
\let\makeshorttitle\maketitlepage
\let\maketitle\maketitlepage

\newwrite\gtoutfile
\long\gdef\makeheadfile{  
{\def\\{, }\def\s{ }
\immediate\openout\gtoutfile head.xxx
\immediate\write\gtoutfile{To: math@arxiv.org}
\immediate\write\gtoutfile{Subject: put OR rep NNNNN:ppppp}
\immediate\write\gtoutfile{--text follows this line--}
\immediate\write\gtoutfile{Proxy-for: \ifx\theasciiauthors\relax
\theauthors\else\theasciiauthors\fi\s<\ifx\theasciiemail\relax\theemail\else\theasciiemail\fi>}
\immediate\write\gtoutfile{\noexpand\\}
\immediate\write\gtoutfile{Authors: \ifx\theasciiauthors\relax
\theauthors\else\theasciiauthors\fi}
{\def\\{ }\immediate\write\gtoutfile{Title: \ifx\theasciititle\relax
\thetitle\else\theasciititle\fi}}
\immediate\write\gtoutfile{Subj-class: GT or SG, GR etc}
\immediate\write\gtoutfile{MSC-class: \theprimaryclass\ifx\thesecondaryclass\relax\else, \thesecondaryclass\fi}
\immediate\write\gtoutfile{Journal-ref: Algebr. Geom. Topol. \thevolumenumber\s
(\thevolumeyear) \startpage-\finishpage}
\immediate\write\gtoutfile{Comments: Published by Algebraic and
Geometric Topology at}
\immediate\write\gtoutfile{\s\s\s  http://www.maths.warwick.ac.uk/agt/AGTVol\thevolumenumber/agt-\thevolumenumber-\thepapernumber.abs.html}
\immediate\write\gtoutfile{\noexpand\\}
\immediate\write\gtoutfile{}
\ifx\theasciiabstract\relax
\immediate\write\gtoutfile{\theabstract}\else
\immediate\write\gtoutfile{\theasciiabstract}\fi
\immediate\write\gtoutfile{}
\immediate\write\gtoutfile{\noexpand\\}
\immediate\write\gtoutfile{}
\immediate\closeout\gtoutfile}}  

\def\maketitlepage{\makeagttitle\makeheadfile}
\let\makeshorttitle\maketitlepage
\let\maketitle\maketitlepage

\def\ifplaintex{\expandafter\ifx\csname documentclass\endcsname\relax}

\def\gtp{{\mathsurround=0pt\it $\cal G\mskip-2mu$eometry \&\ 
$\cal T\!\!$opology $\cal P\!$ublications}}  

\def\Addressesr{\bigskip
{\small \parskip 0pt \leftskip 0pt \rightskip 0pt plus 1fil \def\\{\par}
\sl\theaddress\par
\medskip
\rm Email:\stdspace\tt\theemail\hfill\rm Received:\qua\receiveddate \par}}

\def\recd{{\small Received:\qua\receiveddate\ifx\reviseddate\relax
\else\qquad Revised:\qua\reviseddate\fi\par}} 

\def\lognumber#1{\def\thelognumber{#1}}
\def\volumenumber#1{\def\thevolumenumber{#1}}
\def\volumeyear#1{\def\thevolumeyear{#1}}
\def\papernumber#1{\def\thepapernumber{#1}}
\def\pagenumbers#1#2{\def\startpage{#1}\def\finishpage{#2}}
\def\published#1{\def\publishdate{#1}}

\def\received#1{\def\receiveddate{#1}}

\def\accepted#1{\def\accepteddate{#1}}

\let\\\par\let\thelognumber\relax\let\thevolumenumber\relax
\let\thepapernumber\relax\let\thevolumeyear\relax\let\startpage\relax
\let\finishpage\relax\let\publishdate\relax\let\receiveddate\relax
\let\reviseddate\relax\let\accepteddate\relax\let\theasciititle\relax
\let\theasciiauthors\relax
\let\theasciiabstract\relax

\let\theasciiemail\relax

\ifplaintex
\font\logobig=cmssbx10 scaled 3836
\font\logomed=cmssbx10 scaled 2557
\else
\font\logobig=cmssbx10 scaled 4200
\font\logomed=cmssbx10 scaled 2800
\fi

\long\def\makeagttitle{   
\count0=\startpage
\agt\hfill      
\hbox to 45truept{\vbox to 0pt{\vglue -13truept{\logomed A\kern -.37em{\logobig 
T}\kern -.38em G}\vss}\hss}
\break
{\small Volume \thevolumenumber\ (\thevolumeyear)
\startpage--\finishpage\nl
Published: \publishdate}

\vglue .25truein

{\parskip=0pt\leftskip 0pt plus
1fil\def\\{\par\smallskip}{\Large\bf\thetitle}\par\medskip} \vglue
0.05truein

{\parskip=0pt\leftskip 0pt plus 1fil\def\\{\par}{\sc\theauthors}
\par\medskip}%
 
\vglue 0.03truein 

{\small\leftskip 25truept\rightskip 25truept{\bf Abstract}\stdspace\theabstract

{\bf AMS Classification}\stdspace\theprimaryclass
\ifx\thesecondaryclass\relax\else; \thesecondaryclass\fi\par
{\bf Keywords}\stdspace \thekeywords\par}\vglue 7truept

}   

\ifplaintex
\font\phead=cmsl9 scaled 950
\font\pnum=cmbx10 scaled 913
\font\pfoot=cmsl9 scaled 950
\headline{\vbox to 0pt{\vskip -4.5mm\line{\small\phead\ifnum
\count0=\startpage ISSN 1472-2739 (on-line) 1472-2747 (printed)
\hfill {\pnum\folio}\else\ifodd\count0\def\\{ }%
\ifx\theshorttitle\relax\thetitle\else\theshorttitle\fi\hfill{\pnum\folio}
\else\def\\{ and }{\pnum\folio}\hfill\ifx\theshortauthors\relax\theauthors
\else\theshortauthors\fi\fi\fi}\vss}}
\footline{\vbox to 0pt{\vglue 0mm\line{\small\pfoot\ifnum\count0=\startpage
\copyright\ \gtp\hfill\else
\agt, Volume \thevolumenumber\ (\thevolumeyear)\hfill\fi}\vss}}
\else
\font=cmsl9 scaled 1050
\font\lnum=cmbx10 
\font=cmsl9 scaled 1050
\makeatletter
\def\@oddhead{{\small\lhead\ifnum\count0=\startpage ISSN 1472-2739 
(on-line) 1472-2747 (printed)\hfill {\lnum\number\count0}\else\ifodd\count0
\def\\{ }\ifx\theshorttitle\relax \thetitle \else\theshorttitle\fi\hfill
{\lnum\number\count0}\else\def\\{ and }{\lnum\number\count0}
\hfill\ifx\theshortauthors\relax 
\theauthors\else\theshortauthors\fi\fi\fi}}\def\@evenhead{\@oddhead}
\def\@oddfoot{\small\lfoot\ifnum\count0=\startpage\copyright\ \gtp\hfill\else
\agt, Volume \thevolumenumber\ (\thevolumeyear)\hfill\fi}
\def\@evenfoot{\@oddfoot}
\makeatother
\fi
\let\maketitlepage\makeagttitle
\let\makeshorttitle\maketitlepage
\let\maketitle\maketitlepage

\newwrite\gtoutfile
\long\gdef\makeheadfile{  
{\def\\{, }\def\s{ }
\immediate\openout\gtoutfile head.xxx
\immediate\write\gtoutfile{To: math@arxiv.org}
\immediate\write\gtoutfile{Subject: put OR rep NNNNN:ppppp}
\immediate\write\gtoutfile{--text follows this line--}
\immediate\write\gtoutfile{Proxy-for: \ifx\theasciiauthors\relax
\theauthors\else\theasciiauthors\fi\s<\ifx\theasciiemail\relax\theemail\else\theasciiemail\fi>}
\immediate\write\gtoutfile{\noexpand\\}
\immediate\write\gtoutfile{Authors: \ifx\theasciiauthors\relax
\theauthors\else\theasciiauthors\fi}
{\def\\{ }\immediate\write\gtoutfile{Title: \ifx\theasciititle\relax
\thetitle\else\theasciititle\fi}}
\immediate\write\gtoutfile{Subj-class: GT or SG, GR etc}
\immediate\write\gtoutfile{MSC-class: \theprimaryclass\ifx\thesecondaryclass\relax\else, \thesecondaryclass\fi}
\immediate\write\gtoutfile{Journal-ref: Algebr. Geom. Topol. \thevolumenumber\s
(\thevolumeyear) \startpage-\finishpage}
\immediate\write\gtoutfile{Comments: Published by Algebraic and
Geometric Topology at}
\immediate\write\gtoutfile{\s\s\s  http://www.maths.warwick.ac.uk/agt/AGTVol\thevolumenumber/agt-\thevolumenumber-\thepapernumber.abs.html}
\immediate\write\gtoutfile{\noexpand\\}
\immediate\write\gtoutfile{}
\ifx\theasciiabstract\relax
\immediate\write\gtoutfile{\theabstract}\else
\immediate\write\gtoutfile{\theasciiabstract}\fi
\immediate\write\gtoutfile{}
\immediate\write\gtoutfile{\noexpand\\}
\immediate\write\gtoutfile{}
\immediate\closeout\gtoutfile}}  

\def\maketitlepage{\makeagttitle\makeheadfile}
\let\makeshorttitle\maketitlepage
\let\maketitle\maketitlepage

\def\ifplaintex{\expandafter\ifx\csname documentclass\endcsname\relax}

\def\gtp{{\mathsurround=0pt\it $\cal G\mskip-2mu$eometry \&\ 
$\cal T\!\!$opology $\cal P\!$ublications}}  

\def\Addressesr{\bigskip
{\small \parskip 0pt \leftskip 0pt \rightskip 0pt plus 1fil \def\\{\par}
\sl\theaddress\par
\medskip
\rm Email:\stdspace\tt\theemail\hfill\rm Received:\qua\receiveddate \par}}

\def\recd{{\small Received:\qua\receiveddate\ifx\reviseddate\relax
\else\qquad Revised:\qua\reviseddate\fi\par}} 

\def\lognumber#1{\def\thelognumber{#1}}
\def\volumenumber#1{\def\thevolumenumber{#1}}
\def\volumeyear#1{\def\thevolumeyear{#1}}
\def\papernumber#1{\def\thepapernumber{#1}}
\def\pagenumbers#1#2{\def\startpage{#1}\def\finishpage{#2}}
\def\published#1{\def\publishdate{#1}}

\def\received#1{\def\receiveddate{#1}}

\def\accepted#1{\def\accepteddate{#1}}

\let\\\par\let\thelognumber\relax\let\thevolumenumber\relax
\let\thepapernumber\relax\let\thevolumeyear\relax\let\startpage\relax
\let\finishpage\relax\let\publishdate\relax\let\receiveddate\relax
\let\reviseddate\relax\let\accepteddate\relax\let\theasciititle\relax
\let\theasciiauthors\relax
\let\theasciiabstract\relax

\let\theasciiemail\relax

\ifplaintex
\font\logobig=cmssbx10 scaled 3836
\font\logomed=cmssbx10 scaled 2557
\else
\font\logobig=cmssbx10 scaled 4200
\font\logomed=cmssbx10 scaled 2800
\fi

\long\def\makeagttitle{   
\count0=\startpage
\agt\hfill      
\hbox to 45truept{\vbox to 0pt{\vglue -13truept{\logomed A\kern -.37em{\logobig 
T}\kern -.38em G}\vss}\hss}
\break
{\small Volume \thevolumenumber\ (\thevolumeyear)
\startpage--\finishpage\nl
Published: \publishdate}

\vglue .25truein

{\parskip=0pt\leftskip 0pt plus
1fil\def\\{\par\smallskip}{\Large\bf\thetitle}\par\medskip} \vglue
0.05truein

{\parskip=0pt\leftskip 0pt plus 1fil\def\\{\par}{\sc\theauthors}
\par\medskip}%
 
\vglue 0.03truein 

{\small\leftskip 25truept\rightskip 25truept{\bf Abstract}\stdspace\theabstract

{\bf AMS Classification}\stdspace\theprimaryclass
\ifx\thesecondaryclass\relax\else; \thesecondaryclass\fi\par
{\bf Keywords}\stdspace \thekeywords\par}\vglue 7truept

}   

\ifplaintex
\font\phead=cmsl9 scaled 950
\font\pnum=cmbx10 scaled 913
\font\pfoot=cmsl9 scaled 950
\headline{\vbox to 0pt{\vskip -4.5mm\line{\small\phead\ifnum
\count0=\startpage ISSN 1472-2739 (on-line) 1472-2747 (printed)
\hfill {\pnum\folio}\else\ifodd\count0\def\\{ }%
\ifx\theshorttitle\relax\thetitle\else\theshorttitle\fi\hfill{\pnum\folio}
\else\def\\{ and }{\pnum\folio}\hfill\ifx\theshortauthors\relax\theauthors
\else\theshortauthors\fi\fi\fi}\vss}}
\footline{\vbox to 0pt{\vglue 0mm\line{\small\pfoot\ifnum\count0=\startpage
\copyright\ \gtp\hfill\else
\agt, Volume \thevolumenumber\ (\thevolumeyear)\hfill\fi}\vss}}
\else
\font=cmsl9 scaled 1050
\font\lnum=cmbx10 
\font=cmsl9 scaled 1050
\makeatletter
\def\@oddhead{{\small\lhead\ifnum\count0=\startpage ISSN 1472-2739 
(on-line) 1472-2747 (printed)\hfill {\lnum\number\count0}\else\ifodd\count0
\def\\{ }\ifx\theshorttitle\relax \thetitle \else\theshorttitle\fi\hfill
{\lnum\number\count0}\else\def\\{ and }{\lnum\number\count0}
\hfill\ifx\theshortauthors\relax 
\theauthors\else\theshortauthors\fi\fi\fi}}\def\@evenhead{\@oddhead}
\def\@oddfoot{\small\lfoot\ifnum\count0=\startpage\copyright\ \gtp\hfill\else
\agt, Volume \thevolumenumber\ (\thevolumeyear)\hfill\fi}
\def\@evenfoot{\@oddfoot}
\makeatother
\fi
\let\maketitlepage\makeagttitle
\let\makeshorttitle\maketitlepage
\let\maketitle\maketitlepage

\newwrite\gtoutfile
\long\gdef\makeheadfile{  
{\def\\{, }\def\s{ }
\immediate\openout\gtoutfile head.xxx
\immediate\write\gtoutfile{To: math@arxiv.org}
\immediate\write\gtoutfile{Subject: put OR rep NNNNN:ppppp}
\immediate\write\gtoutfile{--text follows this line--}
\immediate\write\gtoutfile{Proxy-for: \ifx\theasciiauthors\relax
\theauthors\else\theasciiauthors\fi\s<\ifx\theasciiemail\relax\theemail\else\theasciiemail\fi>}
\immediate\write\gtoutfile{\noexpand\\}
\immediate\write\gtoutfile{Authors: \ifx\theasciiauthors\relax
\theauthors\else\theasciiauthors\fi}
{\def\\{ }\immediate\write\gtoutfile{Title: \ifx\theasciititle\relax
\thetitle\else\theasciititle\fi}}
\immediate\write\gtoutfile{Subj-class: GT or SG, GR etc}
\immediate\write\gtoutfile{MSC-class: \theprimaryclass\ifx\thesecondaryclass\relax\else, \thesecondaryclass\fi}
\immediate\write\gtoutfile{Journal-ref: Algebr. Geom. Topol. \thevolumenumber\s
(\thevolumeyear) \startpage-\finishpage}
\immediate\write\gtoutfile{Comments: Published by Algebraic and
Geometric Topology at}
\immediate\write\gtoutfile{\s\s\s  http://www.maths.warwick.ac.uk/agt/AGTVol\thevolumenumber/agt-\thevolumenumber-\thepapernumber.abs.html}
\immediate\write\gtoutfile{\noexpand\\}
\immediate\write\gtoutfile{}
\ifx\theasciiabstract\relax
\immediate\write\gtoutfile{\theabstract}\else
\immediate\write\gtoutfile{\theasciiabstract}\fi
\immediate\write\gtoutfile{}
\immediate\write\gtoutfile{\noexpand\\}
\immediate\write\gtoutfile{}
\immediate\closeout\gtoutfile}}  

\def\maketitlepage{\makeagttitle\makeheadfile}
\let\makeshorttitle\maketitlepage
\let\maketitle\maketitlepage

\def\ifplaintex{\expandafter\ifx\csname documentclass\endcsname\relax}

\def\gtp{{\mathsurround=0pt\it $\cal G\mskip-2mu$eometry \&\ 
$\cal T\!\!$opology $\cal P\!$ublications}}  

\def\Addressesr{\bigskip
{\small \parskip 0pt \leftskip 0pt \rightskip 0pt plus 1fil \def\\{\par}
\sl\theaddress\par
\medskip
\rm Email:\stdspace\tt\theemail\hfill\rm Received:\qua\receiveddate \par}}

\def\recd{{\small Received:\qua\receiveddate\ifx\reviseddate\relax
\else\qquad Revised:\qua\reviseddate\fi\par}} 

\def\lognumber#1{\def\thelognumber{#1}}
\def\volumenumber#1{\def\thevolumenumber{#1}}
\def\volumeyear#1{\def\thevolumeyear{#1}}
\def\papernumber#1{\def\thepapernumber{#1}}
\def\pagenumbers#1#2{\def\startpage{#1}\def\finishpage{#2}}
\def\published#1{\def\publishdate{#1}}

\def\received#1{\def\receiveddate{#1}}

\def\accepted#1{\def\accepteddate{#1}}

\let\\\par\let\thelognumber\relax\let\thevolumenumber\relax
\let\thepapernumber\relax\let\thevolumeyear\relax\let\startpage\relax
\let\finishpage\relax\let\publishdate\relax\let\receiveddate\relax
\let\reviseddate\relax\let\accepteddate\relax\let\theasciititle\relax
\let\theasciiauthors\relax
\let\theasciiabstract\relax

\let\theasciiemail\relax

\ifplaintex
\font\logobig=cmssbx10 scaled 3836
\font\logomed=cmssbx10 scaled 2557
\else
\font\logobig=cmssbx10 scaled 4200
\font\logomed=cmssbx10 scaled 2800
\fi

\long\def\makeagttitle{   
\count0=\startpage
\agt\hfill      
\hbox to 45truept{\vbox to 0pt{\vglue -13truept{\logomed A\kern -.37em{\logobig 
T}\kern -.38em G}\vss}\hss}
\break
{\small Volume \thevolumenumber\ (\thevolumeyear)
\startpage--\finishpage\nl
Published: \publishdate}

\vglue .25truein

{\parskip=0pt\leftskip 0pt plus
1fil\def\\{\par\smallskip}{\Large\bf\thetitle}\par\medskip} \vglue
0.05truein

{\parskip=0pt\leftskip 0pt plus 1fil\def\\{\par}{\sc\theauthors}
\par\medskip}%
 
\vglue 0.03truein 

{\small\leftskip 25truept\rightskip 25truept{\bf Abstract}\stdspace\theabstract

{\bf AMS Classification}\stdspace\theprimaryclass
\ifx\thesecondaryclass\relax\else; \thesecondaryclass\fi\par
{\bf Keywords}\stdspace \thekeywords\par}\vglue 7truept

}   

\ifplaintex
\font\phead=cmsl9 scaled 950
\font\pnum=cmbx10 scaled 913
\font\pfoot=cmsl9 scaled 950
\headline{\vbox to 0pt{\vskip -4.5mm\line{\small\phead\ifnum
\count0=\startpage ISSN 1472-2739 (on-line) 1472-2747 (printed)
\hfill {\pnum\folio}\else\ifodd\count0\def\\{ }%
\ifx\theshorttitle\relax\thetitle\else\theshorttitle\fi\hfill{\pnum\folio}
\else\def\\{ and }{\pnum\folio}\hfill\ifx\theshortauthors\relax\theauthors
\else\theshortauthors\fi\fi\fi}\vss}}
\footline{\vbox to 0pt{\vglue 0mm\line{\small\pfoot\ifnum\count0=\startpage
\copyright\ \gtp\hfill\else
\agt, Volume \thevolumenumber\ (\thevolumeyear)\hfill\fi}\vss}}
\else
\font=cmsl9 scaled 1050
\font\lnum=cmbx10 
\font=cmsl9 scaled 1050
\makeatletter
\def\@oddhead{{\small\lhead\ifnum\count0=\startpage ISSN 1472-2739 
(on-line) 1472-2747 (printed)\hfill {\lnum\number\count0}\else\ifodd\count0
\def\\{ }\ifx\theshorttitle\relax \thetitle \else\theshorttitle\fi\hfill
{\lnum\number\count0}\else\def\\{ and }{\lnum\number\count0}
\hfill\ifx\theshortauthors\relax 
\theauthors\else\theshortauthors\fi\fi\fi}}\def\@evenhead{\@oddhead}
\def\@oddfoot{\small\lfoot\ifnum\count0=\startpage\copyright\ \gtp\hfill\else
\agt, Volume \thevolumenumber\ (\thevolumeyear)\hfill\fi}
\def\@evenfoot{\@oddfoot}
\makeatother
\fi
\let\maketitlepage\makeagttitle
\let\makeshorttitle\maketitlepage
\let\maketitle\maketitlepage

\newwrite\gtoutfile
\long\gdef\makeheadfile{  
{\def\\{, }\def\s{ }
\immediate\openout\gtoutfile head.xxx
\immediate\write\gtoutfile{To: math@arxiv.org}
\immediate\write\gtoutfile{Subject: put OR rep NNNNN:ppppp}
\immediate\write\gtoutfile{--text follows this line--}
\immediate\write\gtoutfile{Proxy-for: \ifx\theasciiauthors\relax
\theauthors\else\theasciiauthors\fi\s<\ifx\theasciiemail\relax\theemail\else\theasciiemail\fi>}
\immediate\write\gtoutfile{\noexpand\\}
\immediate\write\gtoutfile{Authors: \ifx\theasciiauthors\relax
\theauthors\else\theasciiauthors\fi}
{\def\\{ }\immediate\write\gtoutfile{Title: \ifx\theasciititle\relax
\thetitle\else\theasciititle\fi}}
\immediate\write\gtoutfile{Subj-class: GT or SG, GR etc}
\immediate\write\gtoutfile{MSC-class: \theprimaryclass\ifx\thesecondaryclass\relax\else, \thesecondaryclass\fi}
\immediate\write\gtoutfile{Journal-ref: Algebr. Geom. Topol. \thevolumenumber\s
(\thevolumeyear) \startpage-\finishpage}
\immediate\write\gtoutfile{Comments: Published by Algebraic and
Geometric Topology at}
\immediate\write\gtoutfile{\s\s\s  http://www.maths.warwick.ac.uk/agt/AGTVol\thevolumenumber/agt-\thevolumenumber-\thepapernumber.abs.html}
\immediate\write\gtoutfile{\noexpand\\}
\immediate\write\gtoutfile{}
\ifx\theasciiabstract\relax
\immediate\write\gtoutfile{\theabstract}\else
\immediate\write\gtoutfile{\theasciiabstract}\fi
\immediate\write\gtoutfile{}
\immediate\write\gtoutfile{\noexpand\\}
\immediate\write\gtoutfile{}
\immediate\closeout\gtoutfile}}  

\def\maketitlepage{\makeagttitle\makeheadfile}
\let\makeshorttitle\maketitlepage
\let\maketitle\maketitlepage

\lognumber{11}
\volumenumber{1}
\volumeyear{2001}
\papernumber{11}
\pagenumbers{231}{241}
\received{17 December 2000}
\accepted{12 April 2001}
\published{16 April 2001}

\let\oldover\over

\usepackage{amssymb, amsmath, epsf, latexsym}

\let\over\oldover

\def\zz{{\bf Z}}

\newtheorem{thm}{Theorem}[section]

\newtheorem{lemma}[thm]{Lemma}
\newtheorem{prop}[thm]{Proposition}

\def\cfigure#1#2#3{\begin{figure}[ht!] 
\epsfxsize=#1
\centerline{\epsfbox{#2}}
\caption{#3}\end{figure}}

\newenvironment{pf}{\proof}{\endproof}

\begin{document}

\title{Infinite Order Amphicheiral Knots}
\author{Charles Livingston}
\address{Department of Mathematics, Indiana University\\Bloomington, IN 47405,
USA}
\email{livingst@indiana.edu}
\begin{abstract}  
In answer to a question of Long, Flapan constructed an example of a
prime strongly positive amphicheiral knot that is not slice. Long had
proved that all such knots are algebraically slice. Here we show that
the concordance group of algebraically slice knots contains an
infinitely generated free subgroup that is generated by prime strongly
positive amphicheiral knots. A simple corollary of this result is the
existence of positive amphicheiral knots that are not of order two in
concordance.
\end{abstract}
\maketitle


\section{Introduction} 

In 1977 Gordon \cite{go1} asked whether  every class of order two in
the classical knot concordance group can be represented by an
amphicheiral knot.  The question remains open although
counterexamples in higher dimensions are now known to exist
\cite{cm}.  This problem is more naturally stated in terms of {\em
negative} amphicheiral knots, since such knots  represent 2--torsion
in concordance; that is, if $K$ is negative amphicheiral, $ K
\# K$ is slice.
 On the other hand there is no reason to assume that {\em positive}
amphicheiral knots represent 2--torsion.  Surprisingly, until now  no
examples of positive amphicheiral knots that are not 2--torsion  in
concordance have been known.  The first goal of this paper is to
present an example of a   positive amphicheiral knot that is of
infinite order in the concordance group.

Long
\cite{lo} proved that every {\em strongly positive} amphicheiral
knot is algebraically slice.  (This result followed work of Hartley
and Kawauchi \cite{hk} showing a close relationship between strongly
positive amphicheiral knots and slice knots.)   Long also used  an 
example of the author
\cite{li1} to show that such knots need not be slice.   Long's
example was composite and he asked   whether such an example exists
among prime strongly amphicheiral knots.  Flapan provided such an
example in
\cite{fl}.  Here we will prove that examples such as Flapan's are of
infinite order in the concordance group.  This, of course, provides
an example of a positive amphicheiral knot that  is not 2--torsion
in concordance. We will in fact prove much more, describing an
infinite set of strongly positive amphicheiral  knots that are
linearly independent in the concordance group. We restate this as
the main theorem of this paper.

\begin{thm}\label{mainthm} The group of concordance classes of
algebraically slice knots contains an infinitely generated free
subgroup generated by prime strongly positive amphicheiral knots.
\end{thm}

\noindent{\bf Acknowledgement}\qua  This research was motivated by a
question of Alexander Stoimenow  asking the order of Flapan's knot
in the concordance group.
 
\section{Amphicheirality}   We will assume the reader is familiar
with basic knot theory; \cite{rol} serves as a good reference. 
Material concerning Casson--Gordon invariants will be referenced as
needed.  Since the details of orientation and amphicheirality may be
less familiar, we review them briefly here.

For our purposes it is simplest to view a knot as a smooth oriented
pair, $(S, K)$, where
$S$ is  diffeomorphic to $S^3$ and $K$ is diffeomorphic to $S^1$. 
Such a pair is usually abbreviated simply as
$K$.  Associated to $K$ we have three other knots: the {\em mirror}
of $K$, $ K^{*} = (-S ,K)$; the {\em  reverse} of $K$, $K^r = (S,-K)$;
and the (concordance) {\em inverse} of
$K$,
$-K = (-S, -K)$.   

A knot is called  {\em positive amphicheiral} if 
$( S, K)$ is oriented diffeomorphic to
$(-S , K)$, and is called {\em strongly positive amphicheiral} if
that diffeomorphism can be taken to be an involution.  In terms of a
knot diagram a knot is strongly positive amphicheiral if changing
all the crossings of some oriented diagram of the knot results in a
new   diagram that can be transformed back into the original diagram
by a 180 degree rotation.  In Figure 2 we illustrate a strongly
positive amphicheiral knot.

\section{ Casson--Gordon invariant results}\label{cgsection} To
prove that algebraically slice knots are not slice we need
Casson--Gordon invariants. We now summarize the results that we will
use.  The main result of \cite{cg1} is the following.  Let $K$ be a
knot in $S^3$ and let $M_q$ denote its $q$--fold cyclic branched
cover with $q$ a prime power.

\begin{prop} If a knot
$K$ is slice then there is a subgroup
$H \subset H_1(M_q,\zz)$ such that:
\begin{enumerate}
\item  $H^{\perp} = H$ with respect to the linking form on
$H_1(M_q,\zz)$ and, in particular, the order of $H$ is the
square root of the order of $H_1(M_q,\zz)$. 
\item If    $\chi\co H_1(M_q,\zz) \rightarrow \zz_{p^k}$, $p$ prime,
is homomorphism that vanishes on $H$ then the associated
Casson--Gordon invariant satisfies $\sigma(K,\chi) = 0$.  
\item $H$ is invariant under the group of deck transformations.

\end{enumerate}
\end{prop}

The invariance of $H$ is not stated in \cite{cg1} but follows
readily from the construction of $H$. Note that in \cite{cg1} the
invariant $\sigma(K,\chi)$ is denoted
$\sigma_1(\tau(K,\chi))$.  To use this result we do not need the
explicit definition of the Casson--Gordon invariants, but need only
their properties given in  Propositions
\ref{cgsat} and \ref{trivialrep} below.  

We will be constructing   our examples by starting with a knot $K$
and replacing the neighborhood of an unknotted circle $L$ in the
complement of a Seifert surface for
$K$   with the complement of a knot $J$.  The identification map
switches longitude and meridian so that the resulting manifold is
still $S^3$.  The effect of this construction is to tie that portion
of $K$ that passes through $L$ into a knot.  Details can be found,
for example, in \cite{gl1}.  Note that the construction does not
change the Seifert form of $K$. Furthermore, there is a natural
isomorphism between the homology groups of the
$q$--fold cyclic branched covers for any given $q$.  To see this  a
simple Mayer--Vietoris argument can be used:   the branched cover of
the modified knot is built from the branched cover of
$K $ by removing a collection of homology circles (solid tori) and
replacing them with other homology circles (copies of the complement
of $J$). Similarly, there is a correspondence between characters on
the two homology groups.

We denote the new knot by
$K_J$.  In the construction used below, this process is iterated; we
start with a link $L$ and replace each component with a knot
complement.  

In the next proposition
$\sigma_{j/p^k}(J)$ denotes the classical   signature of $J$ given
as the signature of
$(1-\omega)V + (1-\bar{\omega})V^t$, where $V$ is a Seifert matrix
for $J$ and $\omega = e^{ {j 2 \pi i / p^{k} }}$.

\begin{prop}\label{cgsat} Suppose that $K_J$ is obtained from $K$ by
removing the neighborhood of an unknotted circle $L$  and replacing
it with the complement of a knot $J$ as described above.  Then $
\sigma(K_J,\chi) - 
\sigma(K,\chi)  =\sum_{j=0}^{q-1}
\sigma_{\chi(T^j(\tilde{L}))/p^k}(J)$, where the
$\tilde{L}$ is a lift of  $L$, $T$ denotes the generator of the
group of deck transformations, and $\chi \co H_1(M_q, \zz) \rightarrow
\zz_{p^k}$ is a homomorphism.
\end{prop}

\begin{pf}  The proof is contained in \cite{gl1}.  It is based on a
similar result for computing Casson--Gordon invariants of satellite
knots first found by Litherland in
\cite{lit2}
\end{pf} The following result is an immediate consequence of a
result of Litherland in \cite{lit2}.

\begin{prop}\label{trivialrep} If $\chi$ is the trivial character,
then
$\sigma(K,
\chi) = 0$.
\end{prop}
\begin{pf}  In Corollary B2 of \cite{lit2} it is shown that for the
trivial character
$\chi$ one has $\sigma_\xi(K,\chi) = \sum_{\zeta^\nu = \xi} 
\sigma_\zeta(K)  -
\sum_{\omega^\nu = 1}  \sigma_\omega(K) $. Here 
$\sigma_\xi(K,\chi)$ is a more general Casson--Gordon invariant than
what we use here; it is related to  $\sigma(K,\chi)$ and the
proposition follows from the formula  $ \lim_{\xi \to 1}
\sigma_\xi(K,\chi) =
\sigma(K,\chi)$.
\end{pf}

\section{Basic building blocks} Consider the knot $K_{ 2m+1}$
illustrated in Figure 1.  In the figure a Seifert surface
$F$ for
$K_{2m+1}$ is evident.     The bands are twisted in such a way that
the Seifert form is
$$ V_m = \left( \begin{matrix} 0 & m+1 \\
               m & 0 
\end{matrix} \right) . $$  Notice that $K_\mu$ is reversible and
that $-K_\mu = K_{-\mu}$, which has Seifert form
$V_{-m-1}$ since
$-2m -1 = 2(-m-1)+1$.

\cfigure{2.4in}{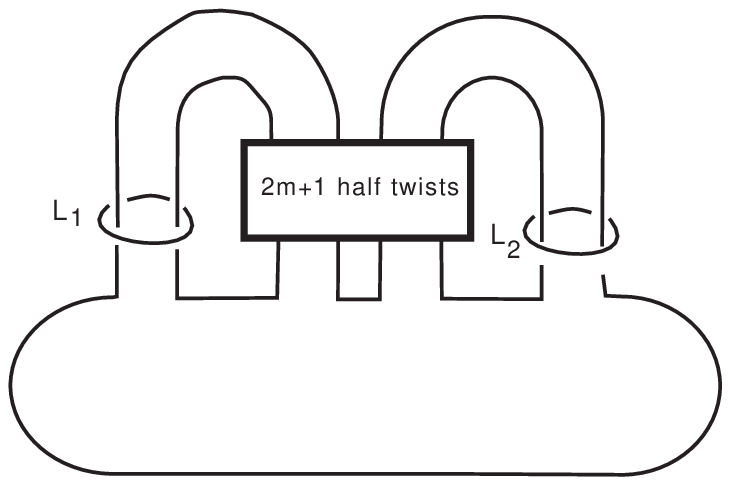}{The knot $K_{2m+1}$}

The link
$L_1
\cup L_2$ will be used later as follows.  A neighborhood of $L_1$
will be removed from
$S^3$ and replaced with the complement of a knot $J $ as described
earlier.  Similarly, a neighborhood of $L_2$ will be removed and
replaced with the complement of $-J$.    This new knot will be
denoted  $K_{J, 2m+1}$.  As mentioned before, a full explication of
this construction is contained in \cite{gl1}.

We will need a detailed understanding of the homology of the 3--fold
cyclic branched cover of
$K_{2m+1}$. Here we give the general result. The proof appears in
detail elsewhere (eg.
\cite{gl1}) but because it is simple we give a fairly complete
outline.

Let $M_q(K_{2m+1})$ denote the $q$--fold cyclic branched cover of
$K_{2m+1}$ and let $\tilde{L}_i$ denote a fixed lift of $L_i$ to
$M_q(K_{2m+1})$, $i=1, 2$.  We will use $\tilde{L}_i$ to denote the
homology classes represented by the $\tilde{L}_i$ also.

\begin{prop}\label{homthm}  $H_1(M_q(K_{2m+1}))= \zz_a \oplus \zz_a$
with $a = |(m+1)^q -  m ^q|$.  The homology is generated by the
classes represented by $\tilde{L}_1$ and
$\tilde{L}_2$.  On homology  the deck transformation $T$ acts by $T(
\tilde{L}_1) = m^*(m+1) \tilde{L}_1$ and     $T( \tilde{L}_2) =
(m+1)^* m   \tilde{L}_2$, where the superscript $(^*)$  denotes the
multiplicative inverse modulo $a$.
 \end{prop}

\begin{pf} A formula of Seifert  (see \cite[Section 8D]{rol}) gives
the  matrix
$\Gamma^q - (\Gamma - I)^q$ as a presentation matrix for the
homology, where $\Gamma = V_m (V_m^t - V_m)$.  That  
$H_1(M_q(K_{2m+1}))= \zz_a \oplus \zz_a$ where $a = |(m+1)^q -  m
^q|$ follows readily. 

A Mayer--Vietoris argument used to compute the homology of the cover
(again see
\cite[Section 8D]{rol}) shows that the   $\tilde{L}_i$ and their
translates generate  the homology of the cover and also gives the
relations
$mT(\tilde{L}_1) = (m+1)(\tilde{L_1})$ and
$(m+1)T(\tilde{L}_2) =  m (\tilde{L_2})$.  Hence, $ T(\tilde{L}_1) =
m^*(m+1)(\tilde{L_1})$ and $ T(\tilde{L}_2) =  (m+1)^*m
(\tilde{L_2})$, as desired. (Notice that $m$ and $m+1$ are
invertible modulo $a$ since $m$ (respectively $m+1$) and
$|(m+1)^q -  m ^q|$ have greatest common divisor 1.)
\end{pf}

\section{A positive  amphicheiral knot of infinite order in
concordance}

Consider the knot $\bar{K}_{2m+1} = K_{2m+1} \# K^*_{2m+1}$
illustrated in Figure 2.  The knot is easily seen to be strongly positive
amphicheiral: reflection through the plane of the page followed by a
180 degree rotation about the axis perpendicular to the plane of the
page gives the desired involution. 

 \cfigure{2.2in}{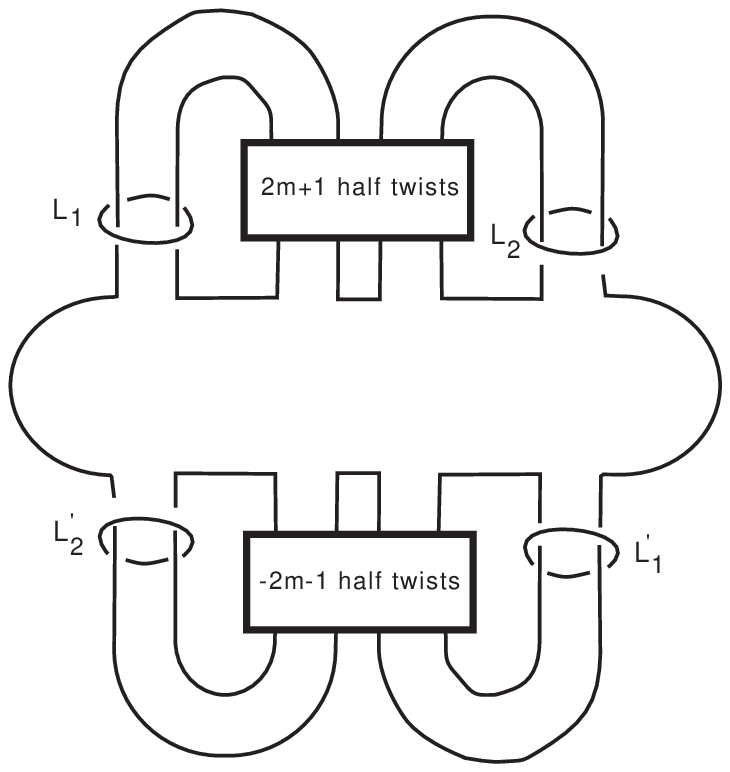}{The strongly positive amphicheiral knot,
$\bar{K}_{2m+1}$}

To form the knot $\bar{K}_{J, 2m+1}$  we replace neighborhoods of
$L_1$ and $L_2'$   with the complement of $J$ and neighborhoods of 
$L_2$ and
$L_1'$   with the complement of $-J$. In this way,  $\bar{K}_{J,
2m+1} = K_{J,2m+1}
\# K^*_{J, 2m+1}$.  (To clarify that the picture is correct, observe
that mirror reflection through the plane of the page, followed by a
180 degree rotation about a point in the center of the figure,
carries the knot complement that replaces $L_1$ to the mirror image
of the knot complement that replaces $L'_1$.  Hence, for this knot
to be strongly positive amphicheiral the knot   complement that
replaces $L_1$ should be the mirror image of the knot that replaces
$L'_1$.)  

Let $p$ be an odd prime dividing $a = (m+1)^3 -  m ^3 $ and suppose
that $p$ has exponent 1 in $(m+1)^3 -  m ^3 $.  That there is an
infinite set of such primes occurring for some $m$ is proved in the
appendix.  In fact, any prime congruent to 1 modulo 3 suffices.

\begin{prop} Suppose that $p$ divides $(m+1)^3 -  m ^3 $ with
exponent $1$. The $p$--torsion in the homology of the $3$--fold
branched cover of
$
\bar{K}_{2m+1}$ is
$(\zz_{p })^4$ generated by the lifts $\{ \tilde{L}_1, \tilde{L}_2,
\tilde{ L}'_1, \tilde{ L}'_2\}$. The deck transformation $T$ has
eigenvalues $\lambda_+ = m^*(m+1)$ with eigenvectors $\{
\tilde{L}_1,  
 \tilde{ L}'_2 \}$ and eigenvalue $\lambda_- = (m+1)^* m$ with
eigenvectors $\{  
\tilde{L}_2,  \tilde{ L}'_1 \}$.
\end{prop}
\begin{pf} This all follows readily from Proposition \ref{homthm}.
\end{pf}

Before proving the next theorem we need to make a simple number
theoretic observation.

\begin{lemma} If $p$ is a prime divisor of $|(m+1)^3 - m^3|$ then
$m^*(m+1) \ne (m+1)^* m$ mod $p$. 
\end{lemma}
\begin{pf} If $m^*(m+1) = (m+1)^* m$ mod $p$ then $m^2 = (m+1)^2$
mod $p$.  It cannot be that $m = m+1$ mod $p$, so we would have
$m = -m - 1$ mod $p$.  In other words, $m = -1/2$ mod $p$. (Since
$p$ is odd, 2 is invertible mod $p$.) Substituting this into
$(m+1)^3 - m^3$ yields $1/4$ mod $p$.  But then $(m+1)^3 - m^3$ is
clearly not divisible by
$p$.
\end{pf}

\begin{thm} Suppose that some odd prime $p$ divides $(m+1)^3 -  m ^3
$ with exponent 1.  For an appropriate choice of $J$ the knot
$\bar{K}_{J, 2m+1}$ has infinite order in the concordance group. 
\end{thm}

\begin{pf} The necessary properties of $J$ will be developed in the
course of the proof.  

Consider the connected sum of
$n $ copies of
$\bar{K}_{J, 2m+1}$. We want to apply the properties of
Casson--Gordon invariants as described in Section \ref{cgsection},
working with the 3--fold cover.

Since the $p$--torsion in $H_1(M_3,\zz)$  is $(\zz_{p })^{4n}$, the
order of $H$ is
$  {p } ^{2n}$. Notice that $(T-\lambda_+ )(T-\lambda_-)$
annihilates the homology of the cover, so
$H$ splits into a summand annihilated by $T - \lambda_+$ and a summand
annihilated  by $T- \lambda_-$.  Here we use that the eigenvalues are
distinct, which follows from the previous lemma.  Call these
$H_+$ and
$H_-$.  One of these has dimension at least 
$  n$; we will assume that it is $H_+$, the other case is identical. 

Any 
$v$ in $H_+$  is a linear combination of the 
$\{ \tilde{L}_1,  
 \tilde{ L}'_2 \}$ associated to each of the $n$ summands.  We want
to see that some such
$v$ involves at least $n$ of these basis elements.  This is an
exercise in elementary linear algebra as follows.  Write out a set
of basis elements for $H_+$ in terms of the full set of
$\{
\tilde{L}_1,  
 \tilde{ L}'_2 \}$.    Writing these as the rows of a matrix (with
at least $n$ rows and exactly $2n$ columns and applying Gauss-Jordan
elimination, combining generators of
$H_+$ and reordering the $\{ \tilde{L}_1,  
 \tilde{ L}'_2 \}$
 yields a generating set with matrix representation
$$  \left( \begin{matrix} 1 & 0 & 0 & \ldots \\
              0  & 1 & 0 & \ldots \\
              0 & 0 & 1 & \ldots \\
              \ldots & & &   
\end{matrix} \right).$$   The sum of these basis elements
corresponding to the rows gives the desired element $v$.

Linking with $v$ defines a character $\chi$ on the homology of the
3--fold cover of
$n\bar{K}_{2m+1}$.  Since $v$ is in the metabolizer, $\chi$ vanishes
on the metabolizer.  Considering just a single building block,
$K_{2m+1}$, we have that the linking form
 satisfies $\mbox{lk}(\tilde{L}_i,\tilde{L}_i) = 0 $, $i = 1$ and
2.   This implies the
$L_1$ and
$L_2$ link nontrivially since the linking form is nonsingular.

From the preceding discussion we see that $\chi$ evaluates
nontrivially on at least $n$ of the $\{\tilde{ L}_2, \tilde{L}'_1
\}$ and trivially on all of the $\{
\tilde{L}_1,  
 \tilde{ L}'_2 \}$.

From Proposition \ref{cgsat} we have that 
$$\sigma(n\bar{K}_{ J,2m+1},\chi) = \sum_i \sum_{j=1}^3  
\sigma_{\alpha_{i,j}/p}(-J)\ \ \ \ +\ \ 
\sigma(n\bar{K}_{2m+1},\chi).$$   In the double summation the index
$i$ runs over a set with at least $n$ terms.  The
$\alpha_{i,j}$ are all nonzero modulo $p$.

The additivity of Casson--Gordon invariants \cite{gi}  applied to
the second summation yields $$\sigma(n\bar{K}_{ J,2m+1},\chi) =
\sum_i \sum_{j=1}^3 \sigma_{\alpha_{i,j}/p}(-J) +
\sum_{i = 1}^n \sigma( \bar{K}_{2m+1},\chi_i).$$ Here the characters
$\chi_i$ are unknown, but notice that the values of $\sigma(
\bar{K}_{2m+1},\chi_i)$ are taken from some finite set of rational
numbers.  Suppose that the maximum value of the absolute values of
this finite set of numbers is $C$.  Then as long as all of the
classical signatures $\sigma_{\alpha /p}( J)$ are greater than
$C/3$, we have that the above sum cannot equal zero.  Finding a knot
$J$ with this property is simple and the proof is complete.
\end{pf}

\noindent{\bf Refinements}\qua  A more detailed analysis of  $\sigma(
\bar{K}_{2m+1},\chi_i)$ can be made to show that in fact these terms
are all 0.  The idea is that the knots $K_{2m+1}$ are doubly slice,
so that the Casson--Gordon invariants must vanish for both
eigenspaces.  Since we do not need this to provide examples for the
main theorem, we do not include the details here.

It is worth noting one particular case that this refinement
addresses, that of Flapan's original example of a nonslice strongly
positive amphicheiral knot in \cite{fl}.  In that case $m = 1$ and
$J$ is the trefoil.  The 3--fold branched  cover has homology $\zz_7
\oplus \zz_7$ and we let $p = 7$.  The associated Casson--Gordon
invariants are then classical 7-signatures of the trefoil.  Although
not all 7--signatures are positive, all are nonnegative and the sum
will necessarily involve some   nonzero terms and is hence
nontrivial.  Hence Flapan's knot is of infinite order in the
concordance group.

\section{Completion of Theorem \ref{mainthm}}

\noindent{\bf An Infinite Family of Linear Independent Examples }

To prove Theorem \ref{mainthm} we need to find an infinite family of
knots such as those constructed in the previous section, all of
which are independent in concordance.  The argument is fairly
simple, similar to the one by Jiang \cite{ji} giving that the
concordance group of algebraically slice knots contains an
infinitely generated free subgroup.

Using the results of the appendix, one can find an infinite sequence
of integers,
$\{m_i\}$, with the property that each $(m_i+1)^3 - m_i^3$ is
divisible by a prime $p_i$ with exponent 1 and no $p_j$ divides
$(m_i+1)^3 - m_i^3$ if $i \ne j$.  

Next construct the $\bar{K}_{J_i,2m_i +1}$ as in the previous
section.  We observe that the $\bar{K}_{J_i,2m_i +1}$ form the desired set, as
follows.  Suppose that some linear combination of these was slice.
Let $\bar{K}_{J_n,2m_n +1}$ be one of the knots in that linear
combination.  Apply the Casson--Gordon result, Proposition
\ref{cgsat}, using characters on the 3--fold cover to $\zz_{p_n}$. 
Using additivity, the value of the Casson--Gordon invariants is
reduced to that on the 3--fold cover of the multiple of
$\bar{K}_{m_n}$.  (The character vanishes on the covers of the other
summands by the choice of the $p_i$ so Proposition \ref{trivialrep}
applies.) The calculations of the previous section show that some
such Casson--Gordon invariants do not vanish.

\medskip
\noindent{\bf Primeness} 

 The proof of Theorem \ref{mainthm} is now completed by showing that
each of the $\bar{K}_{J_i,2m_i +1}$ is concordant to a prime knot. 
In \cite{fl} Flapan showed that one particular example of a
$\bar{K}_{J,2m  +1}$ is concordant to a prime knot that is strongly
positive amphicheiral. Her proof applies in the present setting. 
This completes the proof.

\appendix
\section{Appendix: Prime divisors of $(m+1)^3 - m^3$}
\begin{thm} The set of primes that divide   $F(m) = (m+1)^3 - m^3 =
3m^2 + 3m +1$ for some
$m$ is infinite.
\end{thm}
\begin{pf} The proof is reminiscent of Euclid's proof of the
infinitude of primes. Suppose that the set of such primes is finite,
$\{ {p_i}
\}_{i = 1,
\ldots , n}$.  Let $N = \prod_{i= 1, \ldots , n}p_i$. Consider
$F(N)$. Clearly none of the
$p_i$ divide $F(N) = 3N^2 +3N + 1$, so there is   a prime factor of
$F(N)$ that is not among the
$p_i$, contradicting the definition of $\{ {p_i}
\}$.
\end{pf}

\begin{thm} The set of primes $p$ for which $p$ is a prime divisor
of $F(m)$ with exponent one for some $m$ is infinite.
\end{thm}

\begin{pf} We prove that, in fact, if a prime divides $F(m)$ for
some $m$, then it divides
$F(m + p)$ with exponent 1.  Suppose that $p$ has exponent greater
than 1 in $F(m)$. Use Taylor's theorem to find  $F(m +p) = F(m) +
pF'(m) + p^2F''(m)/2$.  ($F$ is quadratic.)  The first and last
terms are divisible by $p^2$, so $p$ has exponent 1 in $F(m + p)$
unless $p$ divides $F'(m)$.  But, if $F(m)$ and $F'(m)$ have a
common root then
$F(m)$ would have multiple root.  But this can occur if and only if
the discriminant of this quadratic is 0.  In this case the
discriminant is $-3$, so a multiple root does not occur unless $p=3$. 
(As an alternative, to see that $F(m)$ and $F'(m)$ do not have a
common root just note that a simple calculation shows that
$4F(m) - (2m+1)F'(m) = 1$ for all $m$, so that $p$ cannot divide
both $F(m)$ and $F'(m)$.
\end{pf}

Although we don't need a particular description of the primes that
occur, such a calculation is possible.  

\begin{thm} A prime $p$ occurs as a divisor of $F(m)$ with exponent
1 for some $m$ if and only if $p$ is congruent to 1 mod 3.
\end{thm}
\begin{pf} By the previous theorem we need not concern ourselves
with the exponent of
$p$.  Hence, we want the conditions for $3m^2 +3m +1$ to have a root
modulo $p$.  But, using the quadratic formula we quickly see that
this will occur if and only if there is square root of $-3$ in
$\zz_p$.  Using the language of quadratic symbols \cite{la2}, we are
asking for which $p$ we have $\left(   {-3 \over p} \right) = 1$.

From the properties of the quadratic symbol  we have $$\left(   {-3
\over p}
\right) =  \left(   {-1 \over p}
\right) \left(   {3 \over p}
\right)  = (-1)^{{p-1 \over 2}} \left(   {3 \over p}
\right)  .$$ Applying quadratic reciprocity (again, see \cite{la2})
to the second term gives
$$(-1)^{{p-1 \over 2}}
\left(   {3 \over p}
\right) =  (-1)^{{p-1 \over 2}}
\left(   {p \over 3}
\right) (-1)^{({p-1 \over 2})({3-1 \over 2})} = \left(   {p \over 3}
\right).$$ This last term is 1 if and only if $p$ is a square modulo
3, which occurs if and only if
$p$ is congruent to 1 modulo 3.
\end{pf}

\makeatletter\@thebibliography@{CG1}\small\parskip0pt  
plus2pt\relax\makeatother
 
\bibitem[CG1]{cg1}  {\bf A\, Casson}, {\bf  C\, Gordon},
{\it Cobordism of classical knots}, in 
{\it  A la recherche de la Topologie perdue}, ed. by Guillou
and Marin, Progress in Mathematics, Volume 62, 1986, pp.181--199. (Originally
published as Orsay Preprint, 1975.)

\bibitem[CM]{cm} {\bf D\, Coray},  {\bf M\, Daniel},
{\it Knot cobordism and amphicheirality},
 Comment. Math. Helv. 58(4)  (1983)  601--616.

\bibitem[F]{fl} {\bf E\, Flapan},
{\it A prime strongly positive amphicheiral knot which is not
slice},  
 Math. Proc. Cambridge Philos. Soc. 100(3)  (1986)  533--537.

\bibitem[G2]{gi} {\bf P\, Gilmer},
{\it Slice knots in $S\sp{3}$},
  Quart. J. Math. Oxford Ser. (2) 34(135) (1983)
305--322. 

\bibitem[GL1]{gl1}  {\bf P\, Gilmer}, {\bf  C\, Livingston},
{\it  The Casson--Gordon invariant and link concordance}, 
  Topology 31(3)  (1992)  475--492.

\bibitem[G]{go1} {\bf C\, McA\, Gordon}, 
{\it  Problems in Knot Theory},
in {\it Knot Theory}, ed. J. Hausmann, Springer Lecture Notes
no. 685  1977.

\bibitem[HK]{hk} {\bf R\, Hartley},  {\bf  A\, Kawauchi},
{\it Polynomials of amphicheiral knots}, 
  Math. Ann. 243(1) (1979)  63--70.

\bibitem[J]{ji} {\bf B\, Jiang},
{\it  A simple proof that the concordance group of algebraically
slice knots is infinitely generated}, 
  Proc. Amer. Math. Soc. 83 (1981)  189--192.

\bibitem[La]{la2} {\bf S\, Lang},
{\it  Algebraic Number Theory},  
  Addison-Wesley Publishing Co., Reading, Mass., 1968. 

\bibitem[Lt]{lit2} {\bf R\, Litherland}, {\it  Cobordism of satellite knots}, 
in {\it Four--Manifold Theory}, Contemporary Mathematics, eds.
C. Gordon and R. Kirby, American Mathematical Society, Providence RI
1984, 327--362.

\bibitem[L1]{li1} {\bf C\, Livingston},
{\it  Knots which are not concordant to their reverses}, 
 Oxford Q. J. of Math. 34  (1983), 323--328.
 
\bibitem[L2]{li2} {\bf C\, Livingston},
{\it Examples in   concordance}, preprint at
http://front.\-math.\-uc\-davis.\-edu/math.GT/0101035.

\bibitem[Lo]{lo} {\bf D\, D\, Long}, 
{\it Strongly plus-amphicheiral knots are algebraically slice},
  Math. Proc. Cambridge Philos. Soc. 95(2) (1984) 
309--312.

\bibitem[Ro]{rol} {\bf D\, Rolfsen},
  Publish or Perish, Berkeley CA (1976). 

\endthebibliography

\Addressesr
\end{document}